# A Ramanujan enigma involving the first Stieltjes constant

Donal F. Connon

dconnon@btopenworld.com

9 October 2019

**Abstract**

We provide a rigorous formulation of Entry 17(v) in Ramanujan's Notebooks and show how this relates to the first Stieltjes constant.

## 1. Introduction

A curious item appears in Ramanujan's Notebooks ([2] and [3]) which were written by him in India more than 100 years ago.

Entry 17(v) states:

If $0 < x < 1$, then

$$(1.1) \qquad \varphi(x-1) - \varphi(-x) = \pi[\gamma + \log(2\pi)]\cot(\pi x) + 2\pi \sum_{n=1}^{\infty} \sin(2n\pi x)\log n$$

where $\varphi(x)$ is defined as

$$(1.2) \qquad \varphi(x) = \sum_{n=1}^{\infty}\left[\frac{\log n}{n} - \frac{\log(n+x)}{n+x}\right]$$

Ramanujan also notes that

$$(1.3) \qquad \sum_{n=1}^{\infty} \sin(2n\pi x) = \frac{1}{2}\cot(\pi x)$$

In his commentary on Ramanujan's Notebooks, Berndt [2, p.200] has written:

"Of course, Entry 17(v) is meaningless because the series on the right side diverges for $0 < x < 1$. In the midst of his formula, after $\cot(\pi x)$, Ramanujan inserts a parenthetical remark 'for the same limits', the meaning of which we are unable to discern." Unsurprisingly, Berndt goes on to say that (1.3) is also devoid of meaning.

I believe that I may have discovered Ramanujan's modus operandi in Entry 17(v) and can confirm that his two related formulae may be rigorously expressed for $0 < x < 1$ as

$$(1.4) \qquad \varphi(x-1) - \varphi(-x) = \pi[\gamma + \log(2\pi)]\cot(\pi x) + 2\pi \lim_{s \to 1} \sum_{n=1}^{\infty} \frac{\sin(2n\pi x)\log n}{n^{1-s}}$$

and

(1.5) $$\lim_{s \to 1} \sum_{n=1}^{\infty} \frac{\sin(2n\pi x)}{n^{1-s}} = \frac{1}{2}\cot(\pi x)$$

In view of this, we now know the meaning of Ramanujan's parenthetical remark 'for the same limits'.

We shall see later in this paper that

$$\varphi(x-1) - \varphi(-x) = \gamma_1(1-x) - \gamma_1(x)$$

where $\gamma_n(x)$ are the generalised Stieltjes constants defined by (3.1) below.

This results in

(1.6) $$\gamma_1(1-x) - \gamma_1(x) = \pi[\gamma + \log(2\pi)]\cot(\pi x) + 2\pi \lim_{s \to 1} \sum_{n=1}^{\infty} \frac{\sin(2n\pi x)\log n}{n^{1-s}}$$

First of all, we need to set out our building blocks.

## 2. The Hurwitz zeta function $\varsigma(s, x)$

The Hurwitz zeta function $\varsigma(s, x)$ is initially defined for $\operatorname{Re} s > 1$ and $x > 0$ by

(2.1) $$\varsigma(s, x) = \sum_{n=0}^{\infty} \frac{1}{(n+x)^s}$$

Note that $\varsigma(s, x)$ may be analytically continued to the whole $s$ plane except for a simple pole at $s = 1$. For example, Hasse (1898-1979) showed that [13]

(2.2) $$\varsigma(s, x) = \frac{1}{s-1} \sum_{n=0}^{\infty} \frac{1}{n+1} \sum_{k=0}^{n} \binom{n}{k} \frac{(-1)^k}{(k+x)^{s-1}}$$

is a globally convergent series for $\varsigma(s, x)$ and, except for $s = 1$, provides an analytic continuation of $\varsigma(s, x)$ to the entire complex plane.

It may be noted from (2.1) that $\varsigma(s, 1) = \varsigma(s)$.

We easily see from (2.2) that

$$\lim_{s \to 1}[(s-1)\varsigma(s, x)] = \sum_{n=0}^{\infty} \frac{1}{n+1} \sum_{k=0}^{n} \binom{n}{k} (-1)^k$$

and, since $(1-1)^n = \sum_{k=0}^{n} \binom{n}{k}(-1)^k = \delta_{n,0}$, we have

(2.3) $$\lim_{s \to 1}[(s-1)\varsigma(s,x)] = 1$$

which shows that $\varsigma(s,x)$ has a simple pole at $s=1$. This enables us to write the Laurent expansion for the Hurwitz zeta function shown in (3.1) below.

### 3. The generalised Stieltjes constants $\gamma_n(x)$

The generalised Stieltjes constants $\gamma_n(x)$ are the coefficients in the Laurent expansion of the Hurwitz zeta function $\varsigma(s,x)$ about $s=1$

(3.1) $$\varsigma(s,x) = \sum_{n=0}^{\infty} \frac{1}{(n+x)^s} = \frac{1}{s-1} + \sum_{n=0}^{\infty} \frac{(-1)^n}{n!} \gamma_n(x)(s-1)^n$$

We have

(3.2) $$\gamma_0(x) = -\psi(x)$$

where $\psi(x)$ is the digamma function which is the logarithmic derivative of the gamma function $\psi(x) = \frac{d}{du} \log \Gamma(x)$. It is easily seen from the definition of the Hurwitz zeta function that $\varsigma(s,1) = \varsigma(s)$ and accordingly that $\gamma_n(1) = \gamma_n$.

The Stieltjes constants $\gamma_n$ (or the Euler-Mascheroni constants) are the coefficients of the Laurent expansion of the Riemann zeta function $\varsigma(s)$ about $s=1$.

(3.3) $$\varsigma(s) = \frac{1}{s-1} + \sum_{n=0}^{\infty} \frac{(-1)^n}{n!} \gamma_n (s-1)^n$$

Since $\lim_{s \to 1}\left[\varsigma(s) - \frac{1}{s-1}\right] = \gamma$ it is clear that $\gamma_0 = \gamma$. An elementary proof of $\gamma_0(x) = -\psi(x)$ was recently given by the author in [10] (this formula was first obtained by Berndt [1] in 1972).

We note from (3.1) that

(3.4) $$\frac{\partial}{\partial s}[(s-1)\varsigma(s,x)]\bigg|_{s=1} = \gamma_0(x) = -\psi(x)$$

and it is easy to see that

$$\varsigma^{(m)}(s,x) - \varsigma^{(m)}(s,y) = \sum_{n=0}^{\infty} \frac{(-1)^n}{n!}[\gamma_n(x) - \gamma_n(y)]n(n-1)\ldots(n-m+1)(s-1)^{n-m}$$

In the limit $s \to 1$ this becomes

(3.5) $$\lim_{s \to 1}[\varsigma^{(m)}(s,x) - \varsigma^{(m)}(s,y)] = (-1)^m[\gamma_m(x) - \gamma_m(y)]$$

and hence we have for $m \geq 1$

(3.6) $$\gamma_m(x) - \gamma_m(y) = \sum_{n=0}^{\infty}\left[\frac{\log^m(n+x)}{n+x} - \frac{\log^m(n+y)}{n+y}\right]$$

For example, we have

(3.7) $$\gamma_m(x) - \gamma_m = \sum_{n=0}^{\infty}\left[\frac{\log^m(n+x)}{n+x} - \frac{\log^m(n+1)}{n+1}\right]$$

and $m=0$ gives us

(3.8) $$\psi(x) - \gamma = -\sum_{n=0}^{\infty}\left[\frac{1}{n+x} - \frac{1}{n+1}\right]$$

It is easily seen from (2.1) that

(3.9) $$\varsigma(s, 1+x) - \varsigma(s,x) = -\frac{1}{x^s}$$

and (3.5) gives us

(3.10) $$\gamma_m(1+x) - \gamma_m(x) = -\frac{\log^m x}{x}$$

In the particular case $m=0$, we have the familiar formula for the digamma function

(3.11) $$\psi(1+x) = \psi(x) + \frac{1}{x}$$

□

We have for example using (2.1) for $x > 0$ and $1-x > 0$ (i.e. $0 < x < 1$)

(3.12) $$\varsigma(s+1, x) - \varsigma(s+1, 1-x) = \sum_{n=0}^{\infty}\left[\frac{1}{(n+x)^{s+1}} - \frac{1}{(n+1-x)^{s+1}}\right]$$

In the limit as $s \to 0$ we obtain

$$\gamma_0(x) - \gamma_0(1-x) = \sum_{n=0}^{\infty}\left[\frac{1}{n+x} - \frac{1}{n+1-x}\right]$$

and using (3.1) we have

(3.13) $$\varsigma(s+1,x) - \varsigma(s+1,1-x) = \sum_{n=0}^{\infty}\frac{(-1)^n}{n!}[\gamma_n(x) - \gamma_n(1-x)]s^n$$

We see from (2.1) that for $\operatorname{Re}(s) > 1$

(3.14) $$\frac{\partial}{\partial x}\varsigma(s,x) = -s\varsigma(s+1,x)$$

and, by analytic continuation, this holds for all $s$. Using (3.1) we find that

(3.15) $$s\varsigma(s+1,x) = 1 + \sum_{k=0}^{\infty}\frac{(-1)^k}{k!}\gamma_k(x)s^{k+1}$$

Differentiation results in

$$\frac{\partial^{n+1}}{\partial s^{n+1}}[s\varsigma(s+1,x)] = \sum_{k=0}^{\infty}\frac{(-1)^k}{k!}\gamma_k(x)(k+1)k(k-1)\ldots(k-n)s^{k-n}$$

or equivalently

$$\frac{\partial^{n+1}}{\partial s^{n+1}}\frac{\partial}{\partial x}\varsigma(s,x) = \sum_{k=0}^{\infty}\frac{(-1)^k}{k!}\gamma_k(x)(k+1)k(k-1)\ldots(k-n)s^{k-n}$$

We note that the partial derivatives commute in the region where $\varsigma(s,x)$ is analytic and hence we have

$$\frac{\partial^n}{\partial s^n}\frac{\partial}{\partial x}\varsigma(s,x) = \frac{\partial}{\partial x}\frac{\partial^n}{\partial s^n}\varsigma(s,x)$$

Therefore we obtain

(3.16) $$\frac{\partial}{\partial x}\varsigma^{(n+1)}(0,x) = (n+1)(-1)^{n+1}\gamma_n(x)$$

Integration results in

(3.17) $$\int_1^u \gamma_n(x)\,dx = \frac{(-1)^{n+1}}{n+1}\left[\varsigma^{(n+1)}(0,u) - \varsigma^{(n+1)}(0)\right]$$

It is immediately seen that Lerch's formula

(3.18) $$\varsigma'(0,x) = \log \Gamma(x) + \varsigma'(0)$$

arises in the case $n=0$ because $\gamma_0(x) = -\psi(x)$.

We have for $u=2$

(3.19) $$\int_1^2 \gamma_n(x)\,dx = \frac{(-1)^{n+1}}{n+1}\left[\varsigma^{(n+1)}(0,2) - \varsigma^{(n+1)}(0)\right]$$

and, because $\varsigma^{(n)}(0,2) = \varsigma^{(n)}(0)$, we deduce that

(3.20) $$\int_1^2 \gamma_n(x)\,dx = 0$$

and

(3.21) $$\int_0^1 \gamma_n(1+x)\,dx = 0$$

We see from (3.1) that

(3.22) $$\left.\frac{\partial^n}{\partial s^n}[(s-1)\varsigma(s,x)]\right|_{s=1} = (-1)^{n-1} n \gamma_{n-1}(x)$$

## 4. The Ramanujan enigma

**Proposition 4.1**

We have for $0 < x < 1$

(4.1) $$\lim_{s \to 1} \sum_{n=1}^{\infty} \frac{\sin 2n\pi x}{n^{1-s}} = \frac{1}{2}\cot(\pi x)$$

**Proof**

We have the well-known Hurwitz's formula for the Fourier expansion of the Riemann zeta function $\varsigma(s,x)$ as reported in Titchmarsh's treatise [17, p.37]

(4.2) $$\varsigma(s,x) = 2\Gamma(1-s)\left[\sin\left(\frac{\pi s}{2}\right)\sum_{n=1}^{\infty}\frac{\cos 2n\pi x}{(2\pi n)^{1-s}} + \cos\left(\frac{\pi s}{2}\right)\sum_{n=1}^{\infty}\frac{\sin 2n\pi x}{(2\pi n)^{1-s}}\right]$$

where $\mathrm{Re}(s) < 0$ and $0 < x \le 1$. In 2000, Boudjelkha [4] showed that this formula also applies in the region $\mathrm{Re}(s) < 1$. It may be noted that when $x = 1$ this reduces to Riemann's functional equation for $\varsigma(s)$.

Letting $x \to 1-x$ in (4.2) we get for $0 \le x < 1$

(4.3) $$\varsigma(s,1-x) = 2\Gamma(1-s)\left[\sin\left(\frac{\pi s}{2}\right)\sum_{n=1}^{\infty}\frac{\cos 2n\pi x}{(2\pi n)^{1-s}} - \cos\left(\frac{\pi s}{2}\right)\sum_{n=1}^{\infty}\frac{\sin 2n\pi x}{(2\pi n)^{1-s}}\right]$$

and we therefore see that for $0 < x < 1$

(4.4) $$\varsigma(s,x) + \varsigma(s,1-x) = 4\Gamma(1-s)\sin\left(\frac{\pi s}{2}\right)\sum_{n=1}^{\infty}\frac{\cos 2n\pi x}{(2\pi n)^{1-s}}$$

(4.5) $$\varsigma(s,x) - \varsigma(s,1-x) = 4\Gamma(1-s)\cos\left(\frac{\pi s}{2}\right)\sum_{n=1}^{\infty}\frac{\sin 2n\pi x}{(2\pi n)^{1-s}}$$

With regard to (4.5) we note that

$$\lim_{s\to 1}\sum_{n=1}^{\infty}\frac{\sin 2n\pi x}{(2\pi n)^{1-s}} = \lim_{s\to 1}\frac{\varsigma(s,x) - \varsigma(s,1-x)}{4\Gamma(1-s)\cos\left(\frac{\pi s}{2}\right)}$$

and we see that

$$\Gamma(2-s) = \Gamma(1+[1-s])$$
$$= (1-s)\Gamma(1-s)$$

We have

$$\Gamma(1-s)\cos\left(\frac{\pi s}{2}\right) = \Gamma(2-s)\frac{\cos\left(\frac{\pi s}{2}\right)}{1-s}$$

and we first of all consider the limit

$$\lim_{s\to 1}\Gamma(1-s)\cos\left(\frac{\pi s}{2}\right) = \lim_{s\to 1}\frac{\cos\left(\frac{\pi s}{2}\right)}{1-s} = \frac{\pi}{2}$$

where we have employed L'Hôpital's rule. Hence, we have

$$\lim_{s\to 1}\sum_{n=1}^{\infty}\frac{\sin 2n\pi x}{(2\pi n)^{1-s}} = \lim_{s\to 1}\frac{\varsigma(s,x) - \varsigma(s,1-x)}{4\Gamma(1-s)\cos\left(\frac{\pi s}{2}\right)}$$

$$= \frac{\gamma_0(x) - \gamma_0(1-x)}{2\pi}$$

$$= \frac{\psi(1-x) - \psi(x)}{2\pi}$$

where we have used (3.13) in the numerator.

It is well known [16, p.14] that for $0 < x < 1$

$$\psi(1-x)-\psi(x) = \pi\cot(\pi x)$$

and we find that

$$\lim_{s\to 1}\sum_{n=1}^{\infty}\frac{\sin 2n\pi x}{(2\pi n)^{1-s}} = \frac{1}{2}\cot(\pi x)$$

Since

$$\lim_{s\to 1}\sum_{n=1}^{\infty}\frac{\sin 2n\pi x}{(2\pi n)^{1-s}} = \lim_{s\to 1}\sum_{n=1}^{\infty}\frac{\sin 2n\pi x}{n^{1-s}}$$

we obtain a valid interpretation of Ramanujan's formula (1.3)

(4.6) $$\lim_{s\to 1}\sum_{n=1}^{\infty}\frac{\sin 2n\pi x}{n^{1-s}} = \frac{1}{2}\cot(\pi x)$$

It should be specifically noted that

(4.7) $$\lim_{s\to 1}\sum_{n=1}^{\infty}\frac{\sin 2n\pi x}{n^{1-s}} \neq \sum_{n=1}^{\infty}\lim_{s\to 1}\frac{\sin 2n\pi x}{n^{1-s}}$$

because the latter series obviously diverges (since $\lim_{n\to\infty}\sin 2n\pi x \neq 0$). I am aware of at least two interesting mathematical papers where the authors have obtained correct results by boldly starting off with the erroneous assumption that $\sum_{n=1}^{\infty}\sin 2n\pi x = \frac{1}{2}\cot(\pi x)$.

We refer to (4.6) in more detail in the narrative following Proposition 4.3 below.

**Proposition 4.2**

We have for $0 < x < 1$

(4.8) $$\gamma_1(1-x) - \gamma_1(x) = 2\pi\lim_{s\to 1}\sum_{n=1}^{\infty}\frac{\sin(2n\pi x)\log n}{n^{1-s}} + \pi[\gamma + \log(2\pi)]\cot(\pi x)$$

where $\gamma_n(x)$ are the generalised Stieltjes constants defined by (3.1).

**Proof**

We recall (4.4)

$$\varsigma(s,x) - \varsigma(s,1-x) = 4\Gamma(1-s)\cos\left(\frac{\pi s}{2}\right)\sum_{n=1}^{\infty}\frac{\sin 2n\pi x}{(2\pi n)^{1-s}}$$

and substituting the identity

$$\frac{\Gamma(2-s)}{1-s} = \Gamma(1-s)$$

we have the first derivative

(4.9) $\quad \varsigma'(s,x) - \varsigma'(s,1-x) = 4\Gamma(2-s)\frac{\cos\left(\frac{\pi s}{2}\right)}{1-s}\sum_{n=1}^{\infty}\frac{\sin 2n\pi x \log(2\pi n)}{(2\pi n)^{1-s}}$

$$-4\left[\Gamma'(2-s)\frac{\cos\left(\frac{\pi s}{2}\right)}{1-s} + \Gamma(2-s)\frac{d}{ds}\frac{\cos\left(\frac{\pi s}{2}\right)}{1-s}\right]\sum_{n=1}^{\infty}\frac{\sin 2n\pi x}{(2\pi n)^{1-s}}$$

For computational convenience we write $\cos\left(\frac{\pi s}{2}\right) = -\sin\left(\frac{\pi}{2}[1-s]\right)$ so that

$$\frac{d}{ds}\frac{\cos\left(\frac{\pi s}{2}\right)}{1-s} = -\frac{\pi}{2}\frac{d}{ds}\frac{\sin\left(\frac{\pi}{2}[1-s]\right)}{\frac{\pi}{2}[1-s]}$$

Employing the Maclaurin series for the sine function we easily deduce that
$\frac{d}{ds}\frac{\sin\left(\frac{\pi}{2}[1-s]\right)}{\frac{\pi}{2}[1-s]}\bigg|_{s=1} = 0$ and, as noted above, we have $\lim_{s\to 1}\frac{\cos\left(\frac{\pi s}{2}\right)}{1-s} = \frac{\pi}{2}$.

Therefore, taking the limit of (4.9) as $s \to 1$ results in

$$\gamma_1(1-x) - \gamma_1(x) = 2\pi \lim_{s\to 1}\sum_{n=1}^{\infty}\frac{\sin 2n\pi x \log(2\pi n)}{(2\pi n)^{1-s}} - 2\pi\Gamma'(1)\lim_{s\to 1}\sum_{n=1}^{\infty}\frac{\sin 2n\pi x}{(2\pi n)^{1-s}}$$

where we have used (3.5).

Substituting $\Gamma'(1) = -\gamma$ and using (4.1) we obtain

(4.10) $\quad \gamma_1(1-x) - \gamma_1(x) = 2\pi \lim_{s\to 1}\sum_{n=1}^{\infty}\frac{\sin(2n\pi x)\log n}{n^{1-s}} + \pi[\gamma + \log(2\pi)]\cot(\pi x)$

This formula does not provide any new information under the change $x \to 1-x$.

Referring to (3.6) we see that

$$\gamma_1(x) - \gamma_1(y) = \sum_{n=0}^{\infty}\left[\frac{\log(n+x)}{n+x} - \frac{\log(n+y)}{n+1}\right]$$

Ramanujan defined $\varphi(x)$ as

$$\varphi(x) = \sum_{n=1}^{\infty}\left[\frac{\log n}{n} - \frac{\log(n+x)}{n+x}\right] = \sum_{n=0}^{\infty}\left[\frac{\log(n+1)}{n+1} - \frac{\log(n+1+x)}{n+1+x}\right]$$

and thus we have

$$\varphi(x-1)-\varphi(-x) = \sum_{n=0}^{\infty}\left[\frac{\log(n+1-x)}{n+1-x} - \frac{\log(n+x)}{n+x}\right]$$

$$= \gamma_1(1-x) - \gamma_1(x)$$

Hence, we obtain a rigorous formulation of Ramanujan's formula (1.1)

$$\varphi(x-1)-\varphi(-x) = 2\pi \lim_{s\to 1}\sum_{n=1}^{\infty}\frac{\sin(2n\pi x)\log n}{n^{1-s}} + \pi[\gamma+\log(2\pi)]\cot(\pi x)$$

Combining (4.1) and (4.8) we obtain

(4.10.1) $$\gamma_1(1-x)-\gamma_1(x) = 2\pi\lim_{s\to 1}\sum_{n=1}^{\infty}\frac{[\gamma+\log(2\pi n)]\sin(2n\pi x)}{n^{1-s}}$$

□

We show below how (4.10) may be employed to formally derive a known result due to Deninger [11].

Formal integration of (4.10) over $[\frac{1}{2}, u]$, where $0 < u < 1$, results in

(4.11) $$\int_{\frac{1}{2}}^{u}\gamma_1(1-x)\,dx - \int_{\frac{1}{2}}^{u}\gamma_1(x)\,dx = 2\pi\int_{\frac{1}{2}}^{u}\lim_{s\to 1}\sum_{n=1}^{\infty}\frac{\sin(2n\pi x)\log n}{(2\pi n)^{1-s}}\,dx + \pi[\gamma+\log(2\pi)]\int_{\frac{1}{2}}^{u}\cot(\pi x)\,dx$$

Using the elementary integral

$$\int_{\frac{1}{2}}^{u}\gamma_1(x)\,dx = \int_{1}^{u}\gamma_1(x)\,dx - \int_{1}^{\frac{1}{2}}\gamma_1(x)\,dx$$

and

$$\int_{\frac{1}{2}}^{u}\gamma_1(1-x)\,dx = -\int_{\frac{1}{2}}^{1-u}\gamma_1(t)\,dt$$

$$= -\int_{1}^{1-u}\gamma_1(t)\,dt + \int_{1}^{\frac{1}{2}}\gamma_1(t)\,dt$$

we may write the left-hand side of (4.11) as

$$-\int_{1}^{1-u}\gamma_1(x)\,dx - \int_{1}^{u}\gamma_1(x)\,dx + 2\int_{1}^{\frac{1}{2}}\gamma_1(x)\,dx$$

and using (3.17)

$$\int_1^u \gamma_1(x)\,dx = \frac{1}{2}[\varsigma''(0,u) - \varsigma''(0)]$$

this becomes

$$-\frac{1}{2}[\varsigma''(0,u) + \varsigma''(0,1-u)] + \varsigma''(0,\tfrac{1}{2})$$

Using the well-known formula $\varsigma(s,\tfrac{1}{2}) = (2^s - 1)\varsigma(s)$ and Lerch's formula (3.18)

$$\varsigma'(0,x) = \log\Gamma(x) - \frac{1}{2}\log(2\pi)$$

it is readily found that

$$\varsigma''(0,\tfrac{1}{2}) = -\tfrac{1}{2}\log^2 2 - \log(2\pi)\log 2$$

so that the left-hand side of (4.11) may therefore be written as

$$\frac{1}{2}[\varsigma''(0,u) + \varsigma''(0,1-u)] - \frac{1}{2}\log^2 2 - \log(2\pi)\log 2$$

We now assume that formal integration results in

$$\int_{\frac{1}{2}}^u \lim_{s\to 1}\sum_{n=1}^\infty \frac{\sin(2n\pi x)\log n}{(2\pi n)^{1-s}}\,dx = \sum_{n=1}^\infty \int_{\frac{1}{2}}^u \lim_{s\to 1}\frac{\sin(2n\pi x)\log n}{(2\pi n)^{1-s}}\,dx$$

$$= \sum_{n=1}^\infty \lim_{s\to 1}\frac{[(-1)^n - \cos(2n\pi u)]\log n}{(2\pi n)^{2-s}}$$

$$= -\frac{1}{2\pi}\sum_{n=1}^\infty \frac{(-1)^{n+1}\log n}{n} - \frac{1}{2\pi}\sum_{n=1}^\infty \frac{\log n}{n}\cos(2n\pi u)$$

It is well known that [5]

$$\varsigma_a'(1) = -\sum_{n=1}^\infty \frac{(-1)^{n+1}\log n}{n} = \gamma\log 2 - \frac{1}{2}\log^2 2$$

and the right-hand side of (4.11) becomes

$$\gamma\log 2 - \frac{1}{2}\log^2 2 - \sum_{n=1}^\infty \frac{\log n}{n}\cos(2n\pi u) + [\gamma + \log(2\pi)]\log\sin(\pi u)$$

With a little algebra we end up with a known result due to Deninger [11]

(4.12) $$\sum_{n=1}^{\infty} \frac{\log n}{n} \cos(2n\pi u) = -\frac{1}{2}[\varsigma''(0,u) + \varsigma''(0,1-u)] + [\gamma + \log(2\pi)]\log(2\sin(\pi u))$$

which we may write as

(4.12.1) $$\sum_{n=1}^{\infty} \frac{\gamma + \log(2\pi n)}{n} \cos(2n\pi u) = -\frac{1}{2}[\varsigma''(0,u) + \varsigma''(0,1-u)]$$

An elementary integration of (4.12) shows that

(4.13) $$\int_0^1 \varsigma''(0,u)\, du = 0$$

which was noted by Deninger [11].

**Proposition 4.3**

We have for $0 < x < 1$

(4.14) $$\lim_{s \to 1} \sum_{n=1}^{\infty} \frac{\cos 2n\pi x}{n^{1-s}} = -\frac{1}{2}$$

**Proof**

We write (4.4) as

$$\frac{\varsigma(s,x) + \varsigma(s,1-x)}{\Gamma(1-s)} = 4\sin\left(\frac{\pi s}{2}\right) \sum_{n=1}^{\infty} \frac{\cos 2n\pi x}{(2\pi n)^{1-s}}$$

and we have the limit

$$\lim_{s \to 1} \frac{\varsigma(s,x) + \varsigma(s,1-x)}{\Gamma(1-s)} = 4\lim_{s \to 1} \sum_{n=1}^{\infty} \frac{\cos 2n\pi x}{(2\pi n)^{1-s}}$$

Noting that

$$\frac{\varsigma(s,t)}{\Gamma(1-s)} = \frac{(1-s)\varsigma(s,t)}{\Gamma(2-s)}$$

we easily determine the limit

$$\lim_{s \to 1} \frac{\varsigma(s,t)}{\Gamma(1-s)} = -1$$

as noted in [18, p.266].

We therefore obtain

$$\lim_{s \to 1} \sum_{n=1}^{\infty} \frac{\cos 2n\pi x}{(2\pi n)^{1-s}} = -\frac{1}{2}$$

or equivalently

$$\lim_{s \to 1} \sum_{n=1}^{\infty} \frac{\cos 2n\pi x}{n^{1-s}} = -\frac{1}{2}$$

Combining this with (4.1) we obtain

$$\lim_{s \to 1} \sum_{n=1}^{\infty} \left[ \frac{\cos 2n\pi x}{n^{1-s}} + i \frac{\sin 2n\pi x}{n^{1-s}} \right] = \frac{1}{2}[-1 + i \cot(\pi x)]$$

or equivalently

$$\lim_{s \to 1} \sum_{n=1}^{\infty} \frac{e^{2\pi i n x}}{n^{1-s}} = \frac{1}{2}[-1 + i \cot(\pi x)]$$

A good explanation for this interesting (and non-intuitive) result is contained in [14, p.101] and [6]. We have hitherto tacitly assumed that $x$ is a real number. Now let us assume that $x = u + iv$ is a complex number with $v > 0$. This gives us

$$\sum_{n=1}^{\infty} \frac{e^{2\pi i n x}}{n^{1-s}} = \sum_{n=1}^{\infty} \frac{e^{-2\pi n v}}{n^{1-s}} e^{2\pi i n u}$$

We see that this series converges for all $s$ and therefore

$$\lim_{s \to 1} \sum_{n=1}^{\infty} \frac{e^{2\pi i n x}}{n^{1-s}} = \sum_{n=1}^{\infty} e^{-2\pi n v} e^{2\pi i n u}$$

$$= \frac{e^{2\pi i x}}{1 - e^{2\pi i x}}$$

By analytic continuation, this holds true for every $x \in \mathbf{R} - \mathbf{Z}$.

Having regard to (4.6)

$$\lim_{s \to 1} \sum_{n=1}^{\infty} \frac{\sin 2n\pi x}{n^{1-s}} = \frac{1}{2} \cot(\pi x)$$

we note that Nielsen [15, p.80] reports a Fourier series for the cotangent function

(4.15) $$\sum_{n=1}^{\infty} si(2n\pi)\sin(2n\pi x) = \frac{1}{2}\cot(\pi x) - \frac{1}{2}\left[\frac{1}{x} - \frac{1}{1-x}\right]$$

where $si(x)$ is the sine integral function. This could be expressed as

$$\lim_{s \to 1} \sum_{n=1}^{\infty} \frac{si(2n\pi)}{n^{1-s}} \sin(2n\pi x) = \frac{1}{2}\cot(\pi x) - \frac{1}{2}\left[\frac{1}{x} - \frac{1}{1-x}\right]$$

or as

$$\lim_{s \to 1} \sum_{n=1}^{\infty} \frac{si(2n\pi s)}{(2\pi n)^{1-s}} \sin(2n\pi x) = \frac{1}{2}\cot(\pi x) - \frac{1}{2}\left[\frac{1}{x} - \frac{1}{1-x}\right]$$

where we have now also included $s$ in the argument $si(2n\pi s)$.

It might be worthwhile exploring this aspect further.

□

We multiply (4.6) by a Riemann integrable function $f(x)$ and formal integration on $[a,b]$ with $0 < a < b < 1$ results in

$$\int_a^b f(x)\cot(\pi x)\, dx = 2 \int_a^b \lim_{s \to 1} \sum_{n=1}^{\infty} \frac{f(x)\sin 2n\pi x}{n^{1-s}}\, dx$$

This may be compared with [8] and [9] where we showed that

(4.16) $$\int_a^b f(x)\cot(\pi x)\, dx = 2 \sum_{n=1}^{\infty} \int_a^b f(x)\sin(2n\pi x)\, dx$$

for suitably behaved functions.

Similarly, we also showed there that

(4.17) $$\int_a^b f(x)\, dx = 2 \sum_{n=0}^{\infty} \int_a^b f(x)\cos(2n\pi x)\, dx$$

and this corresponds with the integration of (4.14) multiplied by a Riemann integrable function $f(x)$

$$-\frac{1}{2}\int_a^b f(x)\, dx = 2 \int_a^b \lim_{s \to 1} \sum_{n=1}^{\infty} \frac{f(x)\cos 2n\pi x}{n^{1-s}}\, dx$$

One would need to consider the circumstances when

$$\lim_{y \to c} \int_a^b f(x,y)\,dx = \int_a^b \lim_{y \to c} f(x,y)\,dx$$

is valid.

In those two earlier papers [8] and [9] we found the following identity to be extremely useful (which is easily verified by multiplying the numerator and the denominator by the complex conjugate $(1-e^{-ix})$)

$$\frac{1}{1-e^{ix}} = \frac{1}{2} - \frac{i}{2}\frac{\sin x}{1-\cos x} = \frac{1}{2} + \frac{i}{2}\cot(x/2)$$

□

We have

$$\lim_{s \to 1} \frac{\varsigma(s,x) - \varsigma(s,1-x)}{\Gamma(1-s)} = 0$$

Therefore we have

$$\lim_{s \to 1} \cos\left(\frac{\pi s}{2}\right) \sum_{n=1}^{\infty} \frac{\sin 2n\pi x}{(2\pi n)^{1-s}} = 0$$

$$\lim_{s \to 1} \frac{\cos\left(\frac{\pi s}{2}\right)}{1-s}(1-s)\sum_{n=1}^{\infty} \frac{\sin 2n\pi x}{(2\pi n)^{1-s}} = 0$$

$$\lim_{s \to 1}(1-s)\sum_{n=1}^{\infty} \frac{\sin 2n\pi x}{(2\pi n)^{1-s}} = 0$$

which concurs with the obvious limit

$$\lim_{s \to 1}(1-s)\cot(\pi x) = 0$$

The series (1.3) features prominently in the first chapter of Hardy's book, *Divergent Series*, [12] where he demonstrates how it may be used to generate various mathematical formulae.

□

We now differentiate (4.4) with respect to $x$ and obtain

$$\frac{\partial}{\partial x}[\varsigma(s,x) + \varsigma(s,1-x)] = -4\Gamma(1-s)\sin\left(\frac{\pi s}{2}\right)\sum_{n=1}^{\infty} \frac{\sin 2n\pi x}{(2\pi n)^{-s}}$$

Referring to (2.1) we see that

$$\frac{\partial}{\partial x}\varsigma(s,x) = -s\varsigma(s+1,x)$$

$$\frac{\partial}{\partial x}\varsigma(s,1-x) = s\varsigma(s+1,1-x)$$

and we obtain

$$s[\varsigma(s+1,x) - \varsigma(s+1,1-x)] = 4\Gamma(1-s)\sin\left(\frac{\pi s}{2}\right)\sum_{n=1}^{\infty}\frac{\sin 2n\pi x}{(2\pi n)^{-s}}$$

Dividing by $s$ gives us

$$\varsigma(s+1,x) - \varsigma(s+1,1-x) = 4\Gamma(1-s)\frac{\sin\left(\frac{\pi s}{2}\right)}{s}\sum_{n=1}^{\infty}\frac{\sin 2n\pi x}{(2\pi n)^{-s}}$$

Taking the limit $s \to 0$

$$\lim_{s \to 0}[\varsigma(s+1,x) - \varsigma(s+1,1-x)] = 2\pi\lim_{s \to 0}\sum_{n=1}^{\infty}\frac{\sin 2n\pi x}{(2\pi n)^{-s}}$$

we see that this corresponds with (4.1).

## Proposition 4.4

We have for $0 < x < 1$

(4.18) $$2\lim_{s \to 0}\sum_{n=1}^{\infty}\frac{\log(n)\cos 2n\pi x}{(2\pi n)^{-s}} = \psi(x) + \frac{\pi}{2}\cot(\pi x) + [\gamma + \log(2\pi)]$$

**Proof**

We differentiate (4.2) with respect to $s$ to obtain

$$\varsigma'(s,x) = 2\Gamma(1-s)\left[\sin\left(\frac{\pi s}{2}\right)\sum_{n=1}^{\infty}\frac{\log(2\pi n)\cos 2n\pi x}{(2\pi n)^{1-s}} + \cos\left(\frac{\pi s}{2}\right)\sum_{n=1}^{\infty}\frac{\log(2\pi n)\sin 2n\pi x}{(2\pi n)^{1-s}}\right]$$

$$+\pi\Gamma(1-s)\left[\cos\left(\frac{\pi s}{2}\right)\sum_{n=1}^{\infty}\frac{\cos 2n\pi x}{(2\pi n)^{1-s}} - \sin\left(\frac{\pi s}{2}\right)\sum_{n=1}^{\infty}\frac{\sin 2n\pi x}{(2\pi n)^{1-s}}\right]$$

$$-2\Gamma'(1-s)\left[\sin\left(\frac{\pi s}{2}\right)\sum_{n=1}^{\infty}\frac{\cos 2n\pi x}{(2\pi n)^{1-s}} + \cos\left(\frac{\pi s}{2}\right)\sum_{n=1}^{\infty}\frac{\sin 2n\pi x}{(2\pi n)^{1-s}}\right]$$

and then differentiate this with respect to $x$

$$\frac{\partial}{\partial x}\varsigma'(s,x) = 2\Gamma(1-s)\left[-\sin\left(\frac{\pi s}{2}\right)\sum_{n=1}^{\infty}\frac{\log(2\pi n)\sin 2n\pi x}{(2\pi n)^{-s}} + \cos\left(\frac{\pi s}{2}\right)\sum_{n=1}^{\infty}\frac{\log(2\pi n)\cos 2n\pi x}{(2\pi n)^{-s}}\right]$$

$$+\pi\Gamma(1-s)\left[-\cos\left(\frac{\pi s}{2}\right)\sum_{n=1}^{\infty}\frac{\sin 2n\pi x}{(2\pi n)^{-s}}-\sin\left(\frac{\pi s}{2}\right)\sum_{n=1}^{\infty}\frac{\cos 2n\pi x}{(2\pi n)^{-s}}\right]$$

$$-2\Gamma'(1-s)\left[-\sin\left(\frac{\pi s}{2}\right)\sum_{n=1}^{\infty}\frac{\sin 2n\pi x}{(2\pi n)^{-s}}+\cos\left(\frac{\pi s}{2}\right)\sum_{n=1}^{\infty}\frac{\cos 2n\pi x}{(2\pi n)^{-s}}\right]$$

Letting $s \to 0$ gives us

$$\frac{\partial}{\partial x}\varsigma'(0,x) = 2\lim_{s\to 0}\sum_{n=1}^{\infty}\frac{\log(2\pi n)\cos 2n\pi x}{(2\pi n)^{-s}} - \pi\lim_{s\to 0}\sum_{n=1}^{\infty}\frac{\sin 2n\pi x}{(2\pi n)^{-s}} - 2\Gamma'(1)\lim_{s\to 0}\sum_{n=1}^{\infty}\frac{\cos 2n\pi x}{(2\pi n)^{-s}}$$

and hence

$$\frac{\partial}{\partial x}\varsigma'(0,x) = 2\lim_{s\to 0}\sum_{n=1}^{\infty}\frac{\log(2\pi n)\cos 2n\pi x}{(2\pi n)^{-s}} - \frac{\pi}{2}\cot(\pi x) - \gamma$$

where we have employed (4.1) and (4.14)

$$\lim_{s\to 1}\sum_{n=1}^{\infty}\left[\frac{\cos 2n\pi x}{(2\pi n)^{1-s}} + i\frac{\sin 2n\pi x}{(2\pi n)^{1-s}}\right] = \frac{1}{2}[-1 + i\cot(\pi x)]$$

Using Lerch's formula (3.18)

$$\varsigma'(0,x) = \log\Gamma(x) + \varsigma'(0)$$

we see that

$$\frac{\partial}{\partial x}\varsigma'(0,x) = \psi(x)$$

We then obtain

(4.19) $$\psi(x) = 2\lim_{s\to 0}\sum_{n=1}^{\infty}\frac{\log(2\pi n)\cos 2n\pi x}{(2\pi n)^{-s}} - \frac{\pi}{2}\cot(\pi x) - \gamma$$

which may be written as

$$\psi(x) = 2\lim_{s\to 0}\sum_{n=1}^{\infty}\frac{\log(n)\cos 2n\pi x}{(2\pi n)^{-s}} - \frac{\pi}{2}\cot(\pi x) - [\gamma + \log(2\pi)]$$

Letting $x \to 1-x$ we obtain

$$\psi(1-x) = 2\lim_{s\to 0}\sum_{n=1}^{\infty}\frac{\log(2\pi n)\cos 2n\pi x}{(2\pi n)^{-s}} + \frac{\pi}{2}\cot(\pi x) - \gamma$$

and subtraction results in the familiar identity

$$\psi(1-x) - \psi(x) = \pi \cot(\pi x)$$

We also see that

(4.20) $$\psi(x) + \psi(1-x) = -2\gamma + 4\lim_{s \to 0} \sum_{n=1}^{\infty} \frac{\log(2\pi n) \cos 2n\pi x}{(2\pi n)^{-s}}$$

□

Formal integration of (4.20) over $[\tfrac{1}{2}, x]$ results in

$$\log \Gamma(x) - \log \Gamma(1-x) = -2\gamma\left(x - \tfrac{1}{2}\right) + 4\lim_{s \to 0} \sum_{n=1}^{\infty} \frac{\log(2\pi n) \sin 2n\pi x}{(2\pi n)^{1-s}}$$

or

$$\log \Gamma(x) - \log \Gamma(1-x) = -2\gamma\left(x - \tfrac{1}{2}\right) + \frac{2}{\pi} \sum_{n=1}^{\infty} \frac{\log(2\pi n) \sin 2n\pi x}{n}$$

and this corresponds with the formula reported by Deninger [8] in 1984 (which is effectively Kummer's Fourier series for the log gamma function).

□

Formal integration of (4.18) over $[\tfrac{1}{2}, u]$ results in

$$2 \int_{\tfrac{1}{2}}^{u} \lim_{s \to 0} \sum_{n=1}^{\infty} \frac{\log(n) \cos 2n\pi x}{(2\pi n)^{-s}} dx = \int_{\tfrac{1}{2}}^{u} \psi(x) dx + \frac{\pi}{2} \int_{\tfrac{1}{2}}^{u} \cot(\pi x) dx + \left(u - \frac{1}{2}\right)[\gamma + \log(2\pi)]$$

Assuming that the following operations are valid

$$\int_{\tfrac{1}{2}}^{u} \lim_{s \to 0} \sum_{n=1}^{\infty} \frac{\log(n) \cos 2n\pi x}{(2\pi n)^{-s}} dx = \lim_{s \to 0} \int_{\tfrac{1}{2}}^{u} \sum_{n=1}^{\infty} \frac{\log(n) \cos 2n\pi x}{(2\pi n)^{-s}} dx$$

$$= \lim_{s \to 0} \sum_{n=1}^{\infty} \frac{\log(n)}{(2\pi n)^{-s}} \int_{\tfrac{1}{2}}^{u} \cos 2n\pi x \, dx$$

$$= \sum_{n=1}^{\infty} \frac{\log n}{2\pi n} \sin 2n\pi u$$

we obtain

$$\sum_{n=1}^{\infty} \frac{\log n}{\pi n} \sin 2n\pi u = \log \Gamma(u) - \log \Gamma(\tfrac{1}{2}) + \frac{1}{2} \log \sin(\pi u) + \left(u - \frac{1}{2}\right)[\gamma + \log(2\pi)]$$

□

Chakraborty et al, [7, p.426] mention in passing that

$$\sum_{n=1}^{\infty}\left\{[\gamma+\log(2\pi n)]\cos 2n\pi x+\frac{\pi}{2}\sin 2n\pi x\right\}=0 \ (?)$$

but it should be noted that this does <u>not</u> agree with (4.19) which may be written as

$$\psi(x)=2\lim_{s\to 0}\sum_{n=1}^{\infty}\frac{[\gamma+\log(2\pi n)]\cos 2n\pi x+\frac{\pi}{2}\sin 2n\pi x}{(2\pi n)^{-s}}$$

**Proposition 4.5**

We have for $0<x<1$

(4.21) $$\lim_{s\to 1}\sum_{n=1}^{\infty}(-1)^{n+1}\frac{\sin(n\pi x)}{n^{1-s}}=\frac{1}{2}\tan(\pi x/2)$$

(4.22) $$\lim_{s\to 1}\sum_{n=1}^{\infty}(-1)^{n+1}\frac{\cos(n\pi x)}{n^{1-s}}=\frac{1}{2}$$

(4.23) $$\lim_{s\to 1}\sum_{n=0}^{\infty}\frac{\sin(2n+1)\pi x}{(2n+1)^{1-s}}=\frac{1}{\sin\pi x}$$

**Proof**

We showed above that for $0<x<1$

$$\lim_{s\to 1}\sum_{n=1}^{\infty}\frac{\sin(2n\pi x)}{n^{1-s}}=\frac{1}{2}\cot(\pi x)$$

and it is clear that we also have the limits

$$\lim_{s\to 1}\sum_{n=1}^{\infty}\frac{\sin(2n\pi x)}{(2n)^{1-s}}=\frac{1}{2}\cot(\pi x)$$

and

$$\lim_{s\to 1}\sum_{n=1}^{\infty}\frac{\sin(n\pi x)}{n^{1-s}}=\frac{1}{2}\cot(\pi x/2)$$

We have the elementary trigonometric identity

$$\cot(x/2)-\tan(x/2)=2\cot x$$

and, writing this as

$$\frac{1}{2}\cot(x/2)-\cot x=\frac{1}{2}\tan(x/2)$$

shows that

$$\lim_{s \to 1} \sum_{n=1}^{\infty} \left[ \frac{\sin(n\pi x)}{n^{1-s}} - \frac{2\sin(2n\pi x)}{(2n)^{1-s}} \right] = \frac{1}{2} \tan(\pi x / 2)$$

We consider the finite sum

$$\sum_{n=1}^{2N} (-1)^{n+1} a_n = \sum_{n=1}^{2N} a_n - 2\sum_{n=1}^{N} a_{2n}$$

which gives rise to the limit

$$\sum_{n=1}^{\infty} (-1)^{n+1} a_n = \sum_{n=1}^{\infty} [a_n - 2a_{2n}]$$

in the case of suitably convergent series.

Hence, we obtain

$$\lim_{s \to 1} \sum_{n=1}^{\infty} (-1)^{n+1} \frac{\sin(n\pi x)}{n^{1-s}} = \frac{1}{2} \tan(\pi x / 2)$$

and we note that this has a very similar structure to the following Euler summation formula derived by Candelpergher [5a, p.125] in 2017

$$\sum_{n \geq 1}^{E} (-1)^{n+1} \sin(n\pi x) = \frac{1}{2} \tan(\pi x / 2)$$

We also showed that for $0 < x < 1$

$$\lim_{s \to 1} \sum_{n=1}^{\infty} \frac{\cos 2n\pi x}{n^{1-s}} = -\frac{1}{2}$$

and, following the same analysis as above, we find that

$$\lim_{s \to 1} \sum_{n=1}^{\infty} (-1)^{n+1} \frac{\cos(n\pi x)}{n^{1-s}} = \frac{1}{2}$$

which also has a very similar structure to the following Euler summation formula derived by Candelpergher [5a, p.125] in 2017

$$\sum_{n \geq 1}^{E} (-1)^{n+1} \cos(n\pi x) = \frac{1}{2}$$

Using the elementary trigonometric identity

$$\cot(x/2) + \tan(x/2) = \frac{2}{\sin x}$$

it is easily deduced that

$$\lim_{s \to 1} \sum_{n=0}^{\infty} \frac{\sin(2n+1)\pi x}{(2n+1)^{1-s}} = \frac{1}{\sin \pi x}$$

Letting $x = \frac{1}{2}$ in (4.18)

$$2 \lim_{s \to 0} \sum_{n=1}^{\infty} \frac{\log(n) \cos 2n\pi x}{(2\pi n)^{-s}} = \psi(x) + \frac{\pi}{2} \cot(\pi x) + [\gamma + \log(2\pi)]$$

results in

$$\lim_{s \to 0} \sum_{n=1}^{\infty} (-1)^n \frac{\log n}{(2\pi n)^{-s}} = \frac{1}{2} \log\left(\frac{\pi}{2}\right)$$

which corresponds with equation (4.38) in Candelpergher's paper [5a, p.125]

$$\sum_{n \geq 1}^{E} (-1)^n \log n = \frac{1}{2} \log\left(\frac{\pi}{2}\right)$$

## 5. Open access to our own work

This paper contains references to various other papers and, rather surprisingly, most of them are currently freely available on the internet. Surely now is the time that <u>all</u> of <u>our</u> work should be freely accessible by <u>all</u>. The mathematics community should lead the way on this by publishing <u>everything</u> on arXiv, or in an equivalent open access repository. We think it, we write it, so why hide it? You know it makes sense.

Wessex House,
Devizes Road,
Upavon,
Wiltshire SN9 6DL